\documentclass[12pt,a4paper]{article}
\usepackage[russian]{babel}
\let\epsilon\varepsilon
\let\tilde\widetilde
\usepackage{graphicx}
\usepackage{amsfonts}
\begin{document}

\begin{center}
\textbf{Malnormality and centers in one-relator relative presentations.}
\\

\rm
Denis E. Lurye
\\

Faculty of Mechanics and Mathematics
\\

Moscow State University
\\

Moscow 119992, Leninskie Gory, MSU, Russian Federation
\\

\textit{doom1990@mail.ru}
\end{center}

{\large\textbf{1 Introduction.}}
\\

Let $G$ be a group and let a group $\tilde{G}$ be obtained
from the group $G$ by adding one generator and one unimodular relator.

\textit{Definition 1.1.} The word $w=\prod\limits_{i=1}^n g_i
t^{\epsilon_i}\in G\ast \langle t\rangle_\infty$, where
$\epsilon_i=\pm1$, is called \textit{unimodular}, if
$\sum\epsilon_i =1$.

$$\hat{G}=\langle G,t|w=1\rangle:= G\ast \langle t\rangle_\infty/\langle \langle
w \rangle \rangle$$

It is known that if $G$ is a torsion-free group
a significant part of one-relator group theory extends to
such unimodular one-relator relative presentations.
For example:

- $G$ embeds (naturally) into $\hat{G}$ \cite{Kl-27};

- $\hat{G}$ is torsion-free \cite{Fo-Ru};

- $\hat{G}$ is not simple if it does not coincide with $G$ \cite{Kl-hyp};

- $\hat{G}$ almost always (with some known exceptions) contains
a non-abelian free subgroup \cite{Kl-5};

One-relator groups with nontrivial center were investigated particularly \cite{Mu,BaTa,Pi}.
Center of every such a group is an infinite cyclic subgroup (with exception,
when the whole group is a free abelian group of rank 2). In the case of unimodular
one-relator relative presentations $\hat{G}$ almost always has a trivial center
(with some known exceptions)\cite{Kl-09}.

Now let $\tilde{G}$ be a group obtained
from the group $G$ by adding one generator and one relator which is a
proper power of a unimodular word.

$$\tilde{G}=\langle G,t|w^k=1\rangle:= G\ast \langle t\rangle_\infty/\langle \langle
w^k \rangle \rangle$$

In paper \cite{Le} Le proved that if $G$ is a torsion-free group then
$G$ embeds (naturally) into $\tilde{G}$ ($G$ is a
subgroup of $\tilde{G}$) and $\tilde{G}$ is relatively hyperbolic
with respect to $G$.

In paper \cite{Kl-Lu} it is proved that $\tilde{G}$  is relatively hyperbolic
with respect to the subgroup $G$ if $G$ is an involution-free or $k \geq 3$.

\textit{Definition 1.2.} A subgroup $K$ of a group $F$ is called
\textit{malnormal} if $K\cap K^f=\{1\}$ $\forall f \in
F\backslash K$.

D. V. Osin \cite{Os} proved that if some group $R$ is relatively hyperbolic
with respect to its subgroup $H$, then $H\cap H^r$ is finite for all
$r\in R\setminus H$. But $H$ is not always malnormal. For example, if
a group $L$ is finite, $Z$ is an infinite cyclic group and $M$ is an arbitrary nontrivial group,
then $R=L\times(M\ast Z)$ is relatively hyperbolic with respect to $H=L\times M$, but
$H$ is not malnormal. Consider $R=(L\times(M\ast Z))\ast Z$ and
$H=L\times M$, where all the groups are the same as above. Then $R$ is also relatively
hyperbolic with respect to $H$, $H$ is not malnormal, and, moreover,  $R$ does not contain a nontrivial
finite normal subgroup.
\footnote{These examples were courteously provided by D. V. Osin}

The aim of this paper is to prove the following

\textbf{Theorem 1.1.} If
a word $w$ is unimodular, $k\geq2$ and $G$ is an arbitrary group, then
the subgroup $G$ is malnormal in the group $\tilde{G}$  defined by presentation (0).

\begin{center}
$\tilde{G}=\langle G,t|w^k =1\rangle:= G\ast \langle
t\rangle_\infty/\langle \langle w^k \rangle \rangle$, $\qquad (0)$
\end{center}

\textbf{Corollary 1.1.} Suppose $G$ is a nontrivial group, then the group $\tilde {G}$ has a trivial center.

\textit{Proof.} From Theorem 1.1. it follows that the center of the group $\tilde{G}$ is contained in
the subgroup $G$ (consider the condition of the commutation with elements from $G$).

The element $t$ is not contained in $G$, because the left
side of the equality $gt=_{\tilde{G}}1$, where $g\in G$ can not be presented as a product of elements,
conjugated to words $w^{\pm k}$. The reason is that every such a product has the exponent-sum of the
letter $t$ that is multiple of $k$.

Hence the element $t$ does not commutate with nontrivial elements
from the subgroup $G$ in $\tilde {G}$, so we get that the center of $\tilde {G}$ is trivial.

Let us remark that there is another proof of Theorem 1.1.
This proof is shorter then the first one, but it is based on two other independent results and does not work
in case $k=2$, when the group $G$ has involutions. In this paper there are two proofs: the independent proof
that does not use additional facts and the short one (Section 7).

This paper is organized as follows. In Section 2 a new
presentation of the group $\tilde{G}$ is obtained. In Section 4 the relation
between a conjugation of two elements in a certain group and reduced Howie diagrams (Section 3) without exterior
faces with two exterior vertices is presented. This is more general technical problem about an abstract presentation.
In Section 5 all the reduced Howie diagrams without exterior
faces with two exterior vertices over presentation from Section 2 are qualified. In Section 6
the proof is finished.

\textbf{Remark.} Notation $G, \tilde{G}, w, k, t$ means the same as above.

The author thanks A. A. Klyachko and D. V. Osin for useful remarks.
\\

{\large\textbf{2 An algebraic lemma.}}
\\

The following lemma is an easy generalisation
of Lemma 2.1 from \cite{Le};
a similar trick with the change of presentation was
used in \cite{Kl-27} and later
in many other works
(see, e.g., \cite{Fo-Ru, Kl-Lu, Kl-5, FeR98}).

\textbf{Lemma 2.1.} If a word $w=g_1t^{\epsilon_1}\dots g_nt^{\epsilon_n}$
is unimodular and cyclically reduced and $n>1$, then the group
$\tilde{G}$ has a relative presentation of the form
\begin{center}
$\tilde{G}=\langle H,t|\{p^t=p^\varphi,p\in
P\backslash\{1\}\},(ct\prod\limits_{i=0}^m(b_i a_i^t))^k=1\rangle,
\qquad (1)$
\end{center}
where $a_i,b_i,c\in H;P,P^\varphi$ are isomorphic subgroups
of the group $H$,
$\varphi:P\rightarrow P^\varphi$ is an isomorphism between them.
In addition

$1)\quad m\geq0$ (i.e. the product in formula $(1)$ is nonempty);

$2)\quad a_i\notin P, b_i\notin P^\varphi$;

$3)\quad \langle P,a_i\rangle=P\ast \langle p_{i1} a_i
p_{i2}\rangle_{n_i}, \langle P^\varphi,b_i\rangle=P^\varphi\ast
\langle q_{i1} b_i q_{i2}\rangle_{m_i}$ in $H$,
where $p_{i1},p_{i2}\in P; q_{i1},q_{i2}\in P^\varphi$, $n_i$ and
$m_i$ are orders of the elements $(p_{i1} a_i p_{i2})$ and $(q_{i1} b_i
q_{i2})$ respectively.

$4)\quad$ the groups $H$, $P$, and $P^\varphi$ are free products of finitely many
isomorphic copies of $G$:
$H=G^{(0)}*\dots*G^{(s)}$, $P=G^{(0)}*\dots*G^{(s-1)}$,
and $P^\varphi=G^{(1)}*\dots*G^{(s)}$, where $s\ge0$
(if $s=0$, the groups $P$ and $P^\varphi$
are trivial) and the isomorphism $\varphi$ is the shift:
$\left(G^{(i)}\right)^\varphi=G^{(i+1)}$.

\textit{Proof.}
First, we show that $\tilde{G}$ has at least one
presentation of the form (1) satisfying condition 4).
Since
$\sum\epsilon_i=1$, the word $w$ can be written in the form
$$
w=\left(\prod g_i^{t^{k_i}}\right)t.
$$
Conjugating, if necessary, $w$ by $t$, we can assume that $k_i\ge0$.
Setting $g^{(i)}=g^{t^i}$ for $g\in G$, $G^{(i)}=G^{t^i}$, $s=\max k_i$,
and
$c=\prod g_i^{(k_i)}$, we see that $\tilde{G}$ has presentation
\begin{center}
$\tilde{G}=\langle G^{(0)}\ast...\ast
G^{(s)},t|\{(g^{(i)})^t=g^{(i+1)},i=0,...,s-1,g\in G\},(ct)^k=1\rangle$
\end{center}
i.e., a presentation of the form (1) (with $m=-1$) satisfying condition 4).

Now, from all presentations of the form (1) satisfying condition 4)
we choose presentations with minimal $s$, and from all these presentations
with minimal $s$ we choose one with minimal $m$. The obtained
presentation (1) is as required.

Indeed, if $m<0$ (i.e., $w=ct$, where $c\in H$), then $s=0$,
because otherwise we might decrease $s$ replacing
all fragments $g^{(s)}$ in the word $c$ by $(g^{(s-1)})^t$. But the
conditions $m<0$ and $s=0$ mean that the initial word $w$ has the form
$w=ct$, where $c\in G$, which contradicts the assumption $n>1$. Thus,
condition 1) holds.

Condition 2) holds because otherwise in presentation (1)
we might replace a fragment $t^{-1}a_it$ with $a_i\in P$
(or a fragment $tb_it^{-1}$ with $b_i\in P^\varphi$) by $a_i^\varphi$
(or by $b_i^{\varphi^{-1}}$, respectively), thereby decreasing $m$ (and not
increasing $s$).

Let us prove condition 3) for the subgroup $\langle P,a_i\rangle$ (for the subgroup $\langle P^\varphi,b_i\rangle$
the proof is similar).

Since  $a_i\notin P$, we have that  the irreducible form of the element $a_i\in H=G^{(0)}\ast...\ast
G^{(s)}$ is  $\alpha g_{s1}...g_{s2}\beta$, where $\alpha,\beta\in
P;g_{s1},g_{s2}\in G^{(s)}$. If $p_{i1}=\alpha^{-1},p_{i2}=\beta^{-1}$, then
$p_{i1}a_ip_{i2}=g_{s1}...g_{s2}$. Obviously, $\langle
P,a_i\rangle=\langle P,p_{i1}a_ip_{i2}\rangle$, by definition, put
$h_i=p_{i1}a_ip_{i2}$.

Consider a natural homomorphism $\psi$:\quad $P\ast\langle
h_i\rangle_{n_i}\rightarrow \langle P,h_i\rangle$. Clearly, it is surjective. It is also injective.
Indeed, let $s=p_1h_i^{k_1}p_2h_i^{k_2}...p_lh_i^{k_l}p_{l+1}$ be an irreducible form of the element
$x$, $x\neq e, x\in P\ast\langle h_i\rangle_{n_i}$, where $l\geq 1,
k_j$ are integers, each $k_j$ is not multiple of $n_i$, and $p_2,...,p_l$ are not identities in $P$.

Each word $h_i$ begins and ends with a nontrivial element from $G^{(s)}$ as a form (may be reducible) in $H=G^{(0)}\ast ...\ast G^{(s)}$,
hence the reductions may appear only in the subwords $h_i^{k_i}$.
These subwords do not reduce completely, in the reduced form they begin and end with a nontrivial element of the group $G^{(s)}$
(it is sufficient to notice that in a free product
every element is conjugate to the cyclically reduced one). Therefore  $\psi(x)\neq e$ in $H$.

Lemma 2.1. is proved.

\textbf{Corollary 2.1.} Let for some $i$ the equality
$a_i^{n_1}p_1...a_i^{n_s}p_s=1$ or $b_i^{n_1}p^{\varphi}_1...b_i^{n_s}p^{\varphi}_s=1$
holds in the group $H$, where $s\geq 1, n_j \in \mathbb{Z}\setminus\{0\}, p_j\in P$ and $p_j\neq 1$ for all $j \in \{1, \dots, s\}$.
If the orders of all nontrivial elements in $G$ are greater than $l$, then there are more than $l$
letters $a_i$ or $a_i^{-1}$ (respectively, $b_i$ or $b_i^{-1}$) standing successively in the word
$a_i^{n_1}...a_i^{n_s}$ (respectively, $b_i^{n_1}...b_i^{n_s}$).

\textit{Proof.} Without loss of generality we can consider the word $a_i^{n_1}p_1...a_i^{n_s}p_s$. Let us write down the element $a_i$ as
$p_{i1}^{-1}p_{i1}a_ip_{i2}p_{i2}^{-1}$, and the element $a_i^{-1}$ as
$p_{i2}p_{i2}^{-1}a_i^{-1}p_{i1}^{-1}p_{i1}$. Since all the elements
from $P$ in the initial word are nontrivial it follows that there is always a nontrivial element from
$P$ between the elements of type $p_{i1}a_ip_{i2}$ and $p_{i2}^{-1}a_i^{-1}p_{i1}^{-1}$ in the rewritten word.
Therefore the reducing of the rewritten word is possible only if there are more than $l$ elements
$p_{i1}a_ip_{i2}$ or $p_{i2}^{-1}a_i^{-1}p_{i1}^{-1}$ standing in succession.
\\

{\large\textbf{3 Howie diagrams.}}
\\

Definitions and assertions from this section were taken from paper \cite{Kl-hyp}.

Throughout this paper, the term <<surface>> means a closed two-dimensional
oriented surface.

\textit{Definition 3.1.} A \textit{map} $M$ on a surface $S$ is a finite set of continuous mappings
$\{\mu_i:D_i\to S\}$, where $D_i$ is a compact oriented two-dimensional
disk, called the $i$th \textit{face} or \textit{cell} of the map; the boundary
of each face $D_i$ is partitioned into finitely many intervals
$e_{ij}\subset D_i$, called the \textit{pre-edges} of the map, by a nonempty
set of points $c_{ij}\in D_i$, called the {\it corners} of the map.
The images of the corners $\mu_i(c_{ij})$ and the pre-edges $\mu_i(e_{ij})$
are called the {\it vertices} and the {\it edges} of the map, respectively.
It is assumed that

1) the restriction of $\mu_i$ to the interior of each
 face $D_i$ is a homeomorphic embedding preserving orientation; the
 restriction of $\mu_i$ to each pre-edge is a homeomorphic embedding;

2) different edges do not intersect;

3) the images of the interiors of different faces do not intersect;

4) $\bigcup\mu_i(D_i)=S$.

Sometimes, we interpret a map $M$ as a continuous mapping
$M\:\coprod D_i\to S$ from a discrete union of disks onto the surface.

The union of all vertices and edges of a map is a graph on the surface,
called the \textit {$1$-skeleton}. \textit{The degree} of a point from $1$-skeleton
is a number of edges incident to this point, if it is a vertex; if this point belongs
to some edge, then its degree is two. In other words, the degree of point $p$ is
equal to $\mid M^{-1}(p)\mid$.

We say that a corner $c$ is a corner at a vertex $v$ if $M(c)=v$.  There
is a natural cyclic order on the set of all corners at a vertex $v$; we
call two corners at $v$ \textit{adjacent} if they are neighboring with
respect to this order.

By abuse of language, we say that a point or a subset of the surface is
contained in a face $D_i$ if it lies in the image of $\mu_i$. Similarly,
we say that a face $D_i$ is contained in some subset $X\subseteq S$ of
the surface $S$ if $M(D_i)\subseteq X$.

We need also the following simple but useful fact,
sometimes called the combinatorial Gauss--Bonnet formula.

\textbf{Assertion 3.1.} \cite{Ger, Pr, Cam} If each corner $c$ of
a map $M$ on a surface $S$
is assigned a number
$\nu(c)$ (called \textit{the weight} or \textit{the value} of the
corner $c$), then
$$
\sum_v K(v)+\sum_D K(D)+\sum_e K(e)=2\chi(S).
$$
Here the summations are over all vertices $v$
and all cells $D$ of the map
and the values
$K(v)$, $K(D)$, and $K(e)$, called the \textit{curvatures}
of the corresponding vertex, cell, and edge,
are defined by the formula:
$$
K(v)\:=2-\sum_c \nu(c),
\qquad
K(D)\:=2-\sum_c (1-\nu(c)),
\qquad
K(e)\:=0,
$$
where the first sum is over all corners at the vertex $v$, and
the second sum is over all corners of the cell $D$.

Suppose that we have a map $M$ on a surface $S$, the corners of the
map are labeled by elements of a group $H$, and the edges are oriented (in
the figures, we draw arrows on the edges) and labeled by elements of a
set $\{t_1,t_2,\dots\}$ disjoint from the group $H$. The label of a corner
or an edge $x$ is denoted by $\lambda(x)$.

The \textit{label of a vertex} $v$ of such a map is defined by the formula
$$
\lambda(v)=\prod_{i=1}^k \lambda(c_i),
$$
where $c_1,\dots,c_k$ are all corners at $v$ listed clockwise.
The label of a vertex is an element of the group $H$ determined up to
conjugacy.

The \textit{label of a face} $D$ is defined by the formula
$$
\lambda(D)=\prod_{i=1}^k
\bigl(\lambda(M(e_i))\bigr)^{\epsilon_i}\lambda(c_i),
$$
where $e_1,\dots,e_k$ and $c_1,\dots,c_k$ are all pre-edges and all
corners of $D$ listed anticlockwise, the endpoints of $e_i$ are
$c_{i-1}$ and $c_i$ (subscripts are
modulo $k$), and $\epsilon_i=\pm1$ depending on whether the homeomorphism
$M:e_i\mathop\to\limits M(e_i)$ preserves or reverses orientation.
Simply speaking, to obtain the label of a face, we should go around its
boundary anticlockwise, writing out the labels of all corners and edges we
meet; the label of an edge traversed against the arrow should be raised to
the power $-1$.

The label of a face is an element of the group $H*F(t_1,t_2,\dots)$ (the
free product of $H$ and the free group with basis $\{t_1,t_2,\dots\}$)
determined up to a cyclic permutation. More precisely, the right-hand side
of our formula for $\lambda(D)$ is called the {\it label of the face $D$
written starting with the pre-edge $e_1$}.

\textit{Definition 3.2.} Such a labeled map is called a \textit{Howie diagram} (or
simply \textit {diagram}) over a relative presentation
$$\langle H,t_1,t_2,\dots\ |\ w_1=1,w_2=1,\dots\rangle, \qquad (2)$$
if

1) some vertices and faces are distinguished and called \textit{exterior};
the remaining vertices and faces are called \textit{interior};

2) the label of each interior face is a cyclic permutation of one of
the words $w_i^{\pm1}$;

3) the label of each interior vertex is the identity element of $H$.

A diagram is said to be \textit{reduced} if it contains no such
edge $e$ that both faces containing $e$ are interior, these faces
are different and the label of one of this face written
starting with the label of $e$ is inverse
to the label of the other face written ending with
the label of $e$;
such a pair of faces with a common edge is
called a \textit {reducible pair}.

The following lemma is an analogue of the van Kampen lemma for relative
presentations.

\textbf{Lemma 3.1.} (Howie Lemma \cite{How}) The natural mapping from a group $H$ to the group with relative
presentation $(2)$ is noninjective if and only if there exists a
spherical diagram over this presentation with no exterior faces and a
single exterior vertex whose label is not 1 in $G$. A
minimal (with respect to the number of faces) such diagram is
reduced. \hfil\break
If this natural mapping is injective, then we have
the equivalence: the image of an element $u\in
H*F(t_1,t_2,\dots)\setminus \{1\}$ is 1 in the group $(2)$ if and only
if there exists a spherical diagram over this presentation without
exterior vertices and with a single exterior face with label $u$. A
minimal (with respect to the number of faces) such diagram is
also reduced.

Diagrams on the sphere with a single exterior face and no exterior
vertices are also called \textit{disk diagrams}, the boundary of the exterior
face of such a diagram is called the \textit{contour} of the diagram.

Let $\varphi:P\to P^\varphi$ be an isomorphism between two subgroups of a group
$H$.  A relative presentation of the form
$$\langle H,t|\{p^t=p^\varphi; p\in P\backslash\{1\}\}, w_1=1, w_2=1,
\dots\rangle \qquad (3)$$
is called a \textit{$\varphi$-presentation}. A diagram over a $\varphi$-presentation $(3)$ is
called \textit{$\varphi$-reduced} if it is reduced and different interior cells
with labels of the form $p^tp^{-\varphi}$, where $p\in P$, have no common
edges.

\textbf{Lemma 3.2.} \cite{Kl-hyp}  A minimal (with respect to the number of faces)
diagram among all spherical diagrams over a given $\varphi$-presentation
without exterior faces and with a single exterior vertex with nontrivial
label is $\varphi$-reduced. If no such diagrams exists, then a minimal
diagram among all disk diagrams with a given label of contour is
$\varphi$-reduced. In other words, the complete $\varphi$-analogue of Lemma 3.1. is valid.

The idea of the proof is shown in Fig.1.
\\

\unitlength 1.00mm \linethickness{0.4pt}
\begin{picture}(127.00,25.00)
\put(9.00,18.66){\makebox(0,0)[cc]{$p_1$}}
\put(9.34,10.66){\makebox(0,0)[cc]{$p_2$}}
\put(45.34,18.66){\makebox(0,0)[cc]{$p_1^{-\varphi}$}}
\put(45.67,10.66){\makebox(0,0)[cc]{$p_2^{-\varphi}$}}
\put(82.66,14.33){\makebox(0,0)[cc]{$p_1p_2$}}
\put(120.67,14.33){\makebox(0,0)[rc]{$(p_1p_2)^{-\varphi}$}}
\put(62.67,14.66){\makebox(0,0)[cc]{$\to$}}
\put(5.00,15.00){\line(1,0){45.00}}
\put(30.00,25.00){\vector(1,0){1.00}}
\put(30.00,15.00){\vector(1,0){1.00}}
\put(30.00,5.00){\vector(1,0){1.00}}
\put(27.50,15.00){\oval(45.00,20.00)[]}
\put(5.02,15.09){\circle*{2.17}} \put(0.00,19.00){\line(5,-4){5.00}}
\put(5.00,15.00){\line(-2,-5){2.80}}
\put(55.00,20.00){\line(-1,-1){5.00}}
\put(50.00,15.00){\line(2,-5){3.60}}
\put(50.00,15.00){\line(5,-3){5.00}}
\put(50.02,15.09){\circle*{2.17}}
\put(102.00,25.00){\vector(1,0){1.00}}
\put(102.00,5.00){\vector(1,0){1.00}}
\put(99.50,15.00){\oval(45.00,20.00)[]}
\put(77.02,15.09){\circle*{2.17}}
\put(72.00,19.00){\line(5,-4){5.00}}
\put(77.00,15.00){\line(-2,-5){2.80}}
\put(127.00,20.00){\line(-1,-1){5.00}}
\put(122.00,15.00){\line(2,-5){3.60}}
\put(122.00,15.00){\line(5,-3){5.00}}
\put(122.02,15.09){\circle*{2.17}}
\end{picture}

\begin{center}
Fig. 1
\end{center}

Suppose that we have a map on a surface all whose edges are oriented
(e.g., a Howie diagram). Such a map has 4 kinds of corners:
$(++)$, $(--)$, $(+-)$, and $(-+)$ (Fig. 2).

\begin{figure}[h]
\begin{center}
\includegraphics[width=0.75\linewidth]{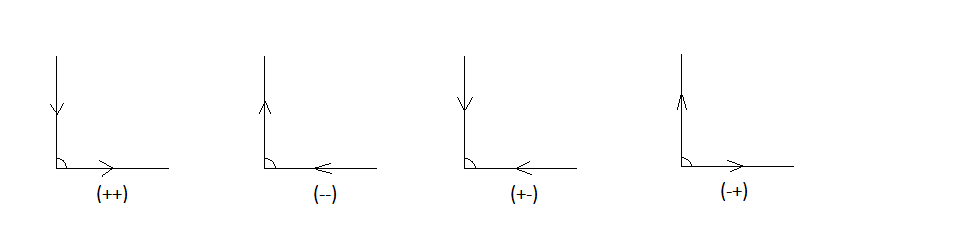}
\end{center}
\begin{center}
Fig. 2
\end{center}
\end{figure}

The following lemma is obvious.

\textbf{Lemma 3.3.} In the anticlockwise listing of the corners at a vertex
$v$, the corners of type $(++)$ alternate with corners of type $(--)$. If
at a vertex $v$ there are no corners of type $(++)$, or, equivalently,
there are no corners of type $(--)$, then either all corners at $v$ are of
type $(+-)$ (in this case, $v$ is called a \textit{sink}), or all corners at
$v$ are of type $(-+)$ (in this case, $v$ is called a \textit{source}).
\\

{\large\textbf{4 The conjugation criterion.}}
\\

\sl
In this section the criterion of conjugation of two elements from the subgroup $H$
in the group $\tilde{H}$ (presentation (3)) is proved. For this purpose
we will define special faces, called 0-cells, and transformations of Howie diagrams.
Further this criterion will be applied to presentation (1).
\\

\rm

Consider a Howie diagram over presentation (3).

\textit{Definition 4.1.} Faces with
labels of the form $p^{-\varphi}p^t$ are called \textit{digons},
the other interior faces are called \textit{large faces}.

\textit{Definition 4.2.} \textit{0-cells}
are the faces shown in Fig. 3, where $x \in H$, 1 is the identity element in $H$ and
the edge on the opposite side of the corner with the identity label has an identity label too.
Such an edge is called a \textit{base edge} of a 0-cell.

\begin{figure}[h]
\begin{center}
\includegraphics[width=0.75\linewidth]{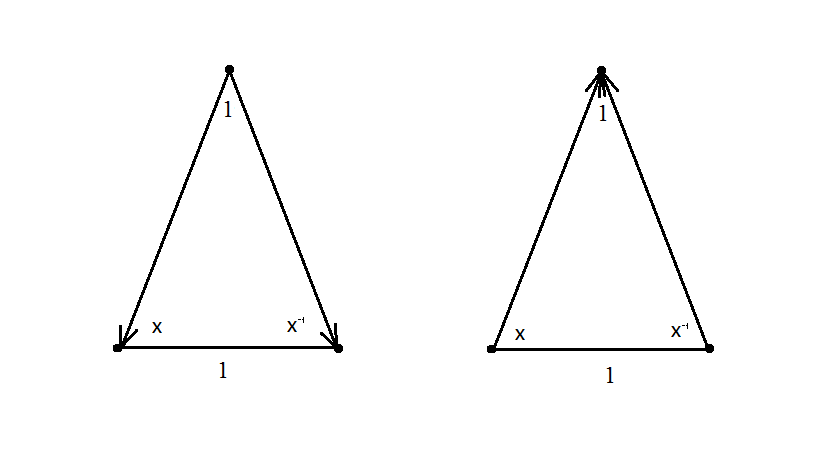}
\end{center}
\begin{center}
Fig. 3
\end{center}
\end{figure}

\textbf{Remark 4.1.} A label of a 0-cell is equal to the identity element in the group $H\ast\langle t\rangle_\infty$.

\textit{Definition 4.3.} Consider a large face, a digon or a 0-cell. Suppose an edge and a corner with identity labels
(\textit{identity edges and corners}) are added in such a face, then the label of the considered face does not change in the free
product $H\ast\langle t\rangle_\infty$. Faces with added identity edges and corners
are called  \textit{extended large faces}, \textit{extended digons}, \textit{extended 0-cells} respectively (Fig. 4).

\begin{figure}[h]
\begin{center}
\includegraphics[width=0.75\linewidth]{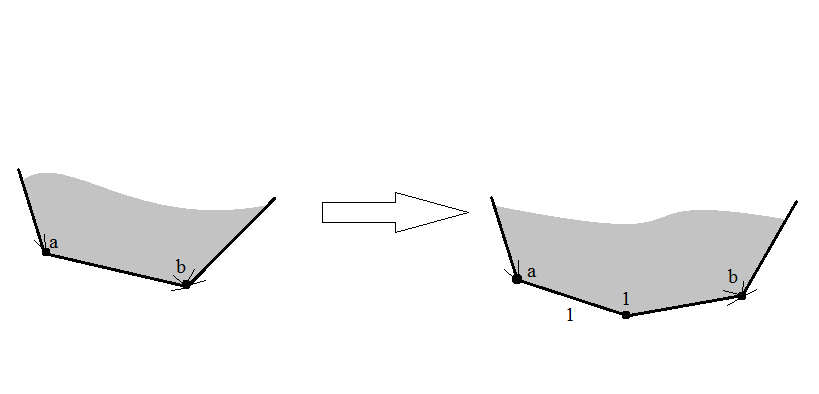}
\end{center}
\begin{center}
Fig. 4
\end{center}
\end{figure}

\textbf{Remark 4.2.} Suppose that the group $H$ is naturally embedded in $\tilde{H}$. In case
of presentation (1) and conditions of Theorem 1.1. it is proved in paper \cite{Kl-Lu}. If
a Howie diagram over the presentation (3) does not have exterior vertices, has the only exterior face and
its interior faces can be extended, then the label of the exterior
face is equal to the identity in $\tilde{H}$. The proof is by induction over the number of
faces in the diagram.

Suppose for some $u \in \tilde{H}$ the equality $u^{-1}hug^{-1}=1$ holds
in the group $H$, where $h,g \in H$.

$H$ is naturally embedded in $\tilde{H}$. Hence from Lemma 3.1. it follows that there is  a $\varphi$-reduced Howie diagram
corresponding to this equality. This diagram is a tree, that consists
of \textit{disks} and \textit{paths}, that connect the disks (Fig. 5: each disk is a diagram, these diagrams
are homeomorphic to 2-dimension disks and consist of interior faces).

\begin{figure}[h]
\begin{center}
\includegraphics[width=0.6\linewidth]{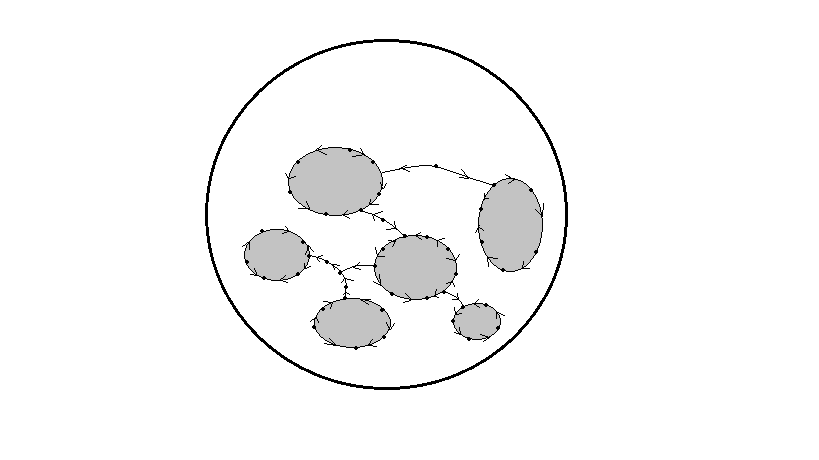}
\end{center}
\begin{center}
Fig. 5
\end{center}
\end{figure}

\textit{Definition 4.4.} Consider a path in a Howie diagram over presentation (3). The product of
all labels of edges and corners (edge and corner labels) in the group $\tilde{H}$ multiplied in a natural order from the beginning to the end of the path is called \textit{a label} of the path. This differs from the definition of a face or a vertex, where a label was an element of the free product.
If a label is considered as an element of the group $H\ast\langle t\rangle_\infty$, then it is called
\textit{a word}. Henceforth the term <<corner>> is more general, it means a set of corners of a diagram at the same vertex
standing successively. By definition all the corner labels are read with a positive power (clockwise) (Fig. 6).

\begin{figure}[h]
\begin{center}
\includegraphics[width=0.6\linewidth]{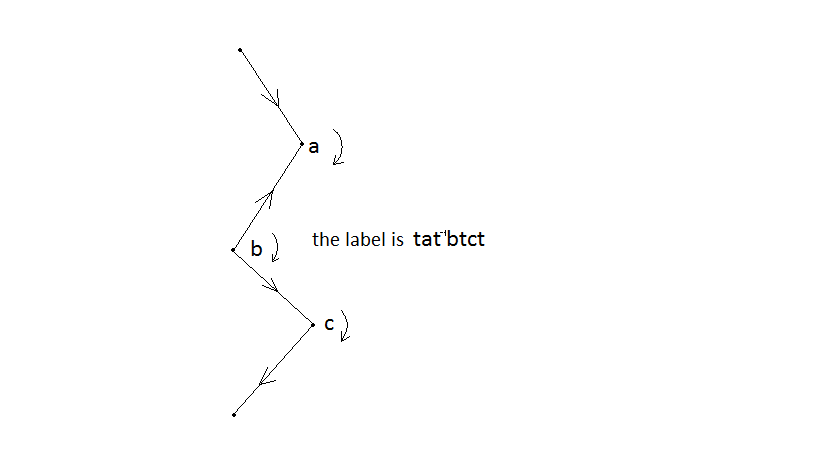}
\end{center}
\begin{center}
Fig. 6
\end{center}
\end{figure}

Now we are going to <<thicken>> paths, i.e. to obtain a diagram that we will be able to
<<glue>> together with itself by the path with label $u$.

\textit{Definition 4.5.} The set of points of the 1-skeleton is called a \textit{marked graph} if all
these points have at least two preimages in the boundary of the exterior face.

Clearly, every marked graph is a finite union of trees.

\textit{Definition 4.6.} A vertex of a marked graph is called \textit{a boundary point}
if this point is contained in a boundary of some interior face. A vertex of a marked graph
that is not a boundary point is called \textit{an interior point} if it has degree greater than two.

\textit{Definition 4.7.} If all the points (boundary and interior) are deleted from the marked
graph, then the obtained parts are called \textit{marked paths}.

Let us replace every interior point by an $n$-gon that has sides with identity labels (identity sides), where $n$ is
the degree of the point. Such an $n$-gon is constructed from $n$ 0-cells with a common vertex.

\textit{Definition 4.8.} Corners at a point of degree $n$ are called \textit{initial corners}.
Let labels of corners at vertices of every $n$-gon be the same as the respective labels of initial corners
(Fig. 7: $n=4$, respective corners are connected). Henceforth edges without arrows have identity labels.
Nontrivial labels at the vertices of an $n$-gon boundary can appear at the corners described above
and at the base edges of 0-cells of an $n$-gon. Let us choose a random 0-cell and begin to assign labels
to the corners at the base edges of 0-cells clockwise in such a way that labels of all vertices of an $n$-gon become
identities. It is possible because the vertex of degree $n$ is interior, i. e. its
label is trivial in the group $H$.

\begin{figure}[h]
\begin{center}
\includegraphics[width=0.6\linewidth]{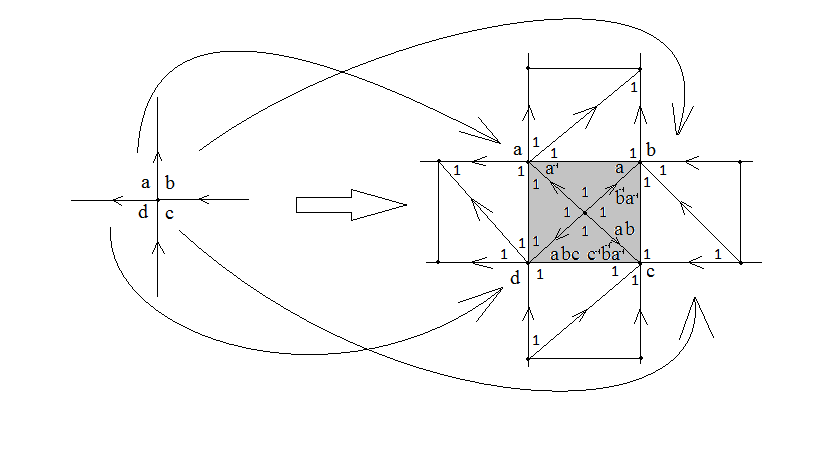}
\end{center}
\begin{center}
Fig. 7
\end{center}
\end{figure}

A boundary point has to be replaced by an identity edge (Fig. 8). This operation
implies existence of identity edges, corners with identity labels and extended faces. If an end of a marked path is an interior
point of degree $n$, then after the <<thikening>> it will become a side of an $n$-gon.

For the <<thikening>> of a marked path we have to duplicate every its edge and connect the respective vertices with the identity edges.
Let us divide each of the obtained quadrangles into two 0-cells. There is a natural order of all these 0-cells generated by the direction
of the path. Let the corner labels at the base edges of odd 0-cells be inverse to the respective corner labels of the path and the corner labels
at the base edges of even 0-cells be trivial (Fig. 8). An exception can be only at the first or the last 0-cell, there all labels can be
nontrivial.

Obviously, the process of <<thikening>> does not add new exterior vertices to the diagram and does not change the label of the exterior face (i. e. does not add identity corners and labels). 

\begin{figure}[h]
\begin{center}
\includegraphics[width=0.6\linewidth]{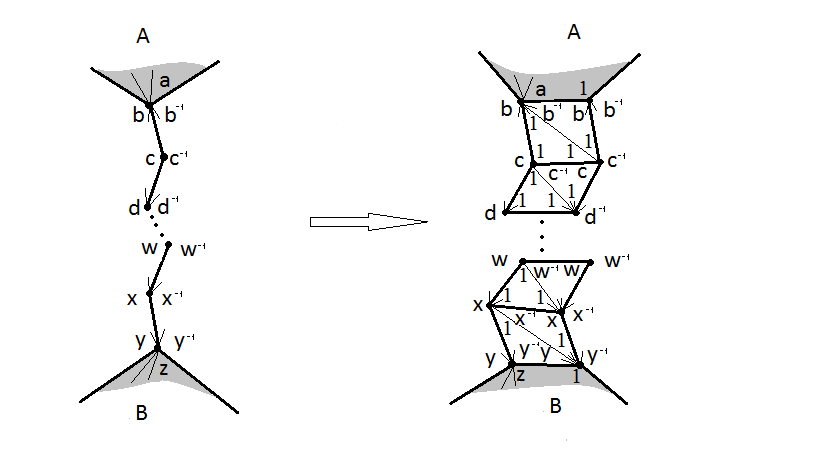}
\end{center}
\begin{center}
Fig. 8
\end{center}
\end{figure}

\textbf{Theorem 4.1.} Elements $h$ and $g$ of the group $H$ are conjugate in the group $\tilde{H}$
($u^{-1}hug^{-1}=1$) $\Leftrightarrow$ these elements are conjugate in $H$ or there are $n \in \mathbb{N}$
$\varphi$-reduced spherical Howie diagrams over presentation (3) without exterior faces with two exterior vertices with labels
$h_0=h$ and $h_1^{-1}$, $h'_1$ and $h_2^{-1}$, ... $h'_{n-1}$ and $h_n^{-1}=g^{-1}$, where $h_i$ and
$h'_i$ are conjugate in $H$ $\forall i \in \{1,2,.., n-1\}$.

\textit{Proof.} Implication $\Leftarrow$ follows from Lemma 3.1. Let us prove another implication.

Consider $u$ as a word from the group $H\ast \langle t\rangle_\infty$. If $u \in H$ then
$h$ and $g$ are conjugate in the group $H$, so we may assume that $u \notin H$.

From Lemma 3.1. and the <<thikening>> process of paths it follows that there is a spherical Howie diagram
over presentation (3) with 0-cells, in this diagram there are two exterior vertices with labels $h$ and $g^{-1}$ and
there are no exterior faces. Further we are going to start a process. Suppose we have $k$ spherical
Howie diagrams over presentation (3) with 0-cells without exterior faces, there are two exterior vertices with
labels $h$ and $h_1^{-1}$, $h'_1$ and $h_2^{-1}$, ...
$h'_{k-1}$ and $g^{-1}$ in every diagram, where $h_i$ and $h'_i$ are conjugate in $H$ $\forall i
\in \{1,2,.., k-1\}$. Further assume that there are no identity edges in some diagram and
there is a reducible pair of large faces in this diagram(a pair of digons can be removed the same way as in Lemma 3.2.).

\textit{Definition 4.9.} \textit{A hole} is a face that is obtained from a reducible pair of large faces
by deleting a path relative to which the labels of large faces are inverse.

\textbf{Remark 4.3.} A label of every hole is equal to the identity element in the free product $H\ast \langle t\rangle_\infty$.

Consider a reducible pair of large faces. Let us replace it by a hole, and fill the hole with 0-cells. We connect corners with inverse labels
using identity edges. We divide each derived quadrangle into two 0-cells. Similarly to the case of the <<thikening>> of paths
there is a natural order of 0-cells. Let the labels of the corners at the base edges of odd 0-cells be the same as at the respective
corners of the hole and the labels of the corners at the base edges of even 0-cells be identity elements (Fig. 9: in this example all
the vertices are interior).

\begin{figure}[h]
\begin{center}
\includegraphics[width=0.6\linewidth]{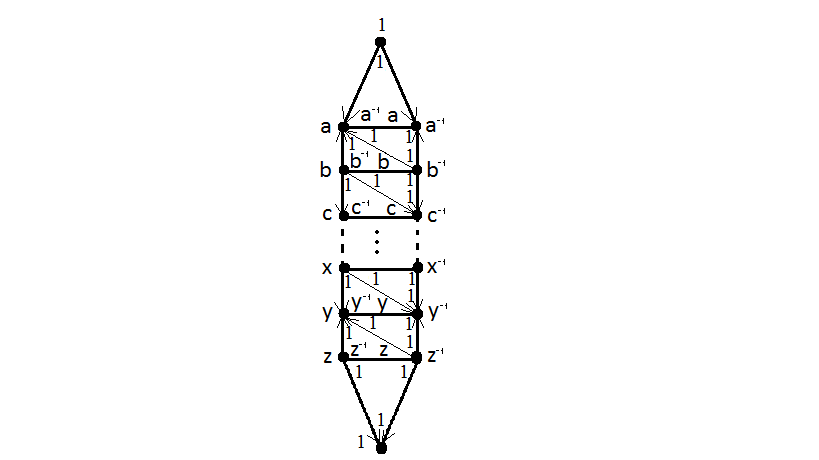}
\end{center}
\begin{center}
Fig. 9
\end{center}
\end{figure}

Now suppose that there are identity edges in the diagram. There are two possibilities.

\begin{enumerate}

\item An identity edge begins and ends with different vertices. Let us pull
this edge into a vertex. If the face that contains the identity edge is a non-extended 0-cell, then
after pulling it will become an edge. If this face is an extended large face, extended digon or extended
0-cell, then after pulling two corners with labels  $a\in H$ and $1$ (the identity in the group $H$)
will be replaced by one corner with the label $a$ (Fig. 10).

\begin{figure}[h]
\begin{center}
\includegraphics[width=0.6\linewidth]{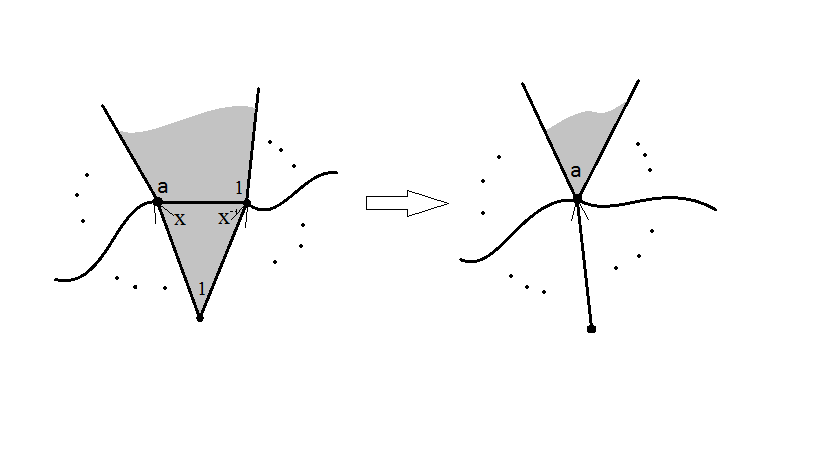}
\end{center}
\begin{center}
Fig. 10
\end{center}
\end{figure}

It is clear that the number of exterior vertices of the diagram does not increase.

\item An identity edge begins and ends with the same vertex. Then the identity edge divides
the diagram into two components. Let $A$ and $B$ be two exterior vertices of the diagram.
Two cases are possible: the exterior vertices are in the same component (Fig.
11a)) or these vertices are in different components (Fig. 11b)). If one of the
exterior vertices coincides with the beginning and the end of the identity edge (Fig. 11c)), then
it is classified as case a).

\begin{figure}[h]
\begin{center}
\includegraphics[width=0.6\linewidth]{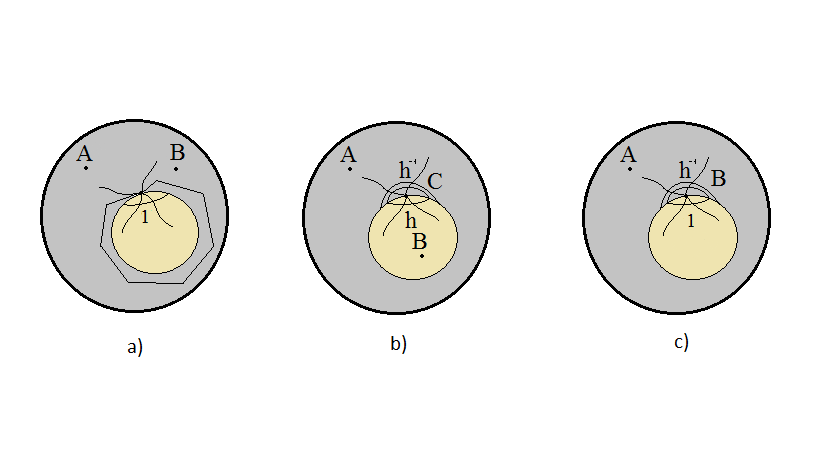}
\end{center}
\begin{center}
Fig. 11
\end{center}
\end{figure}

In case a) the component that contains the exterior vertices is a spherical diagram with two exterior
vertices. Indeed, the product of labels marked by the arc (Fig.
11a)) equals the identity in $H$, it follows from Remark 4.2. Therefore one of the components corresponds to a trivial correlation.

In case b) there are two spherical diagrams, each diagram has two exterior vertices. In the first diagram
these vertices are $A$ and $C$, in the second one these vertices are $C$ and $B$. Note that $C$ is interior in the
initial diagram, hence the labels of two vertices corresponding to $C$ in the initial diagram
satisfy Theorem 4.1. (this means that the label of the considered vertex in the first diagram and the reversed label of
the considered vertex in the second diagram are conjugate in $H$) because these labels are determined up to a cyclic
permutation.
\end{enumerate}

To finish the proof we need to explain the order of operations. Consider one of the spherical diagrams.
If it has an identity edge, then we remove it as we did above, in the other case we search for a reducible pair and remove it
using 0-cells. If the considered diagram is reduced we operate with the next one. The described process is finite because the number
of faces in the initial diagram is finite and therefore the number of identity edges that can appear in this process is finite too.
This implies that finally we will obtain a sequence of reducible Howie diagrams without 0-cells, that will satisfy Theorem 4.1.
\\

{\large\textbf{5 The reduced diagrams case.}}
\\

\sl
In this section all $\varphi$-reduced Howie diagrams without exterior
faces with two exterior vertices over presentation (1) will be qualified.
We use the combinatorial Gauss--Bonnet formula to prove this fact, and its
application will be the same as in paper \cite{Kl-Lu}.
\\

\rm

Consider presentation (1) and a diagram over this presentation.
Faces with a label of the form $p^{-\varphi}p^t$ are called \textit{digons},
the other interior faces are called \textit{large faces}, these terms correspond
to Definition 4.1.

\textbf{Theorem 5.1.} In terms of Theorem 1.1. the only $\varphi$-reduced
Howie diagram over the presentation (1) without exterior faces with two exterior
vertices can be a diagram with only one interior face that is a digon with associated edges (Fig. 9). This
diagram is called \textit{a degenerated digon}.

\begin{figure}[h]
\begin{center}
\includegraphics[width=0.75\linewidth]{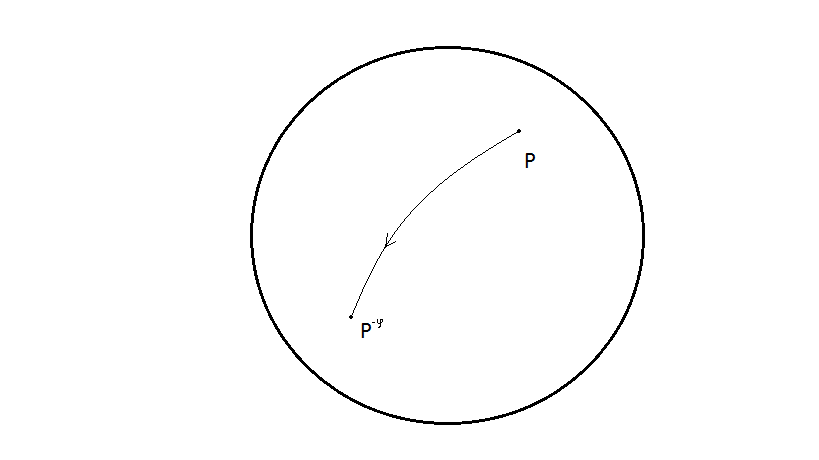}
\end{center}
\begin{center}
Fig. 12
\end{center}
\end{figure}

\textit{Proof.} Let us consider a $\varphi$-reduced Howie diagram without exterior
cells with two exterior vertices over the presentation (1).

A digon is called \textit{special} if
its both neighboring faces are interior and
one of its corners (called \textit{positive})
is adjacent with corners of types $(++)$ and $(--)$
(Fig. 13).
Note that the other corner of a special digon
(called \textit{negative}) is automatically
non-adjacent with corners of type $(++)$ and $(--)$.

\begin{figure}[h]
\begin{center}
\includegraphics[width=0.75\linewidth]{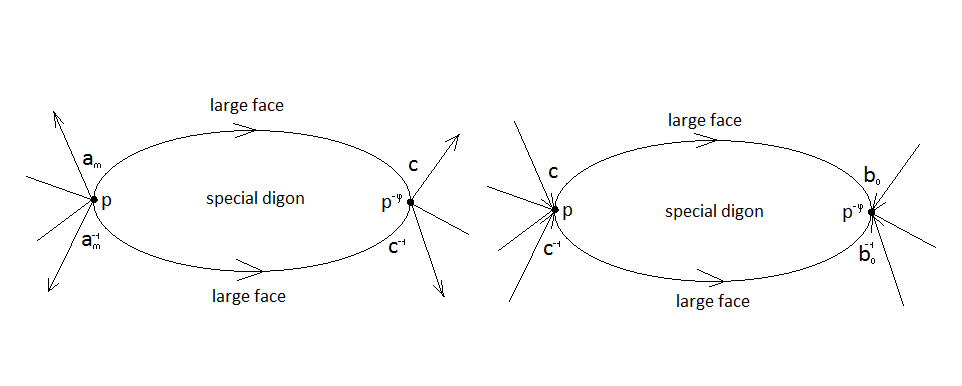}
\end{center}
\begin{center}
Fig. 13
\end{center}
\end{figure}

Let us assign a value (weight)
$\nu(\gamma)$ to each corner $\gamma$ of the diagram by the following
rule:
$$
\nu(\gamma)=\cases{
0 &if $\gamma$ is a corner of a nonspecial digon
\cr
&or a corner of type $(++)$ or $(--)$ of an
interior face (the label of such a corner is $c^{\pm1}$);
\cr
-1 &if $\gamma$ is a negative corner of a special digon;
\cr
1, &otherwise.
\cr
}
$$

Let us calculate the curvatures of vertices and faces according to the
weight test (see Section 3). For faces, we have
$$
K(\hbox{digon})=0,
\quad
K(\hbox{large face})=2-k.
$$
For a vertex $v$, the curvature is
$$
K(v)=2+n-l-{\bf p},
\eqno{(4)}
$$
where $l$ is the number of corners of types $(+-)$ and $(-+)$ of
large faces,
${\bf p}$
is the number of positive corners of special digons,
$n$ is the number of negative corners of special digons
(all corners are at the vertex $v$).

Each negative corner of a special digon is adjacent to
two corners of type $(+-)$ or $(-+)$ of large faces
(by the definition of special digons),
and no corner of type $(+-)$ or $(-+)$ can be adjacent to two
negative corners (since otherwise, the corresponding large face would have both
a corner of type $(++)$ and a corner of type $(--)$).
Therefore, $l\ge 2n$.

Note also that corners of types $(++)$ and $(--)$ at
a non-boundary vertex alternate (Lemma 3.3.) and
cannot be adjacent (since the diagram is reduced):
between two such corners there must be a corner of weight 1
(either a corner of type $(+-)$ or $(-+)$ of a large face or a
positive corner of a special digon).
Taking into account the preceding remark about negative corners,
we conclude that the sum of weights of corners lying between corners of
type $(++)$ and $(--)$ (if we list them clockwise around the vertex $v$)
is at least one (Fig. 14, left). Therefore, a non-boundary
vertex with positive curvature must be either a source or a sink and, for such vertex,
${\bf p}=0$,
and
either $n=1$ and $l=2$ or $n=0$ and $l=1$ or $n=0$ and $l=0$
($n<2$, since otherwise, formula (4) and the inequality
$l\ge 2n$ mentioned above would give a nonpositive curvature).
See Fig. 14, the boldface digits denote the values of corners.

\begin{figure}[h]
\begin{center}
\includegraphics[width=0.75\linewidth]{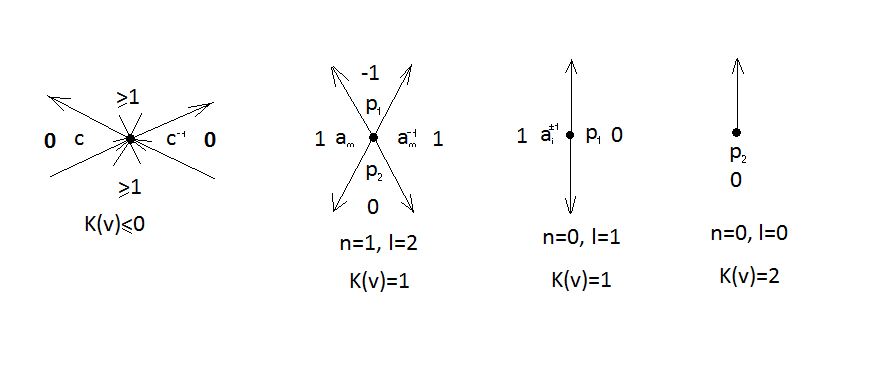}
\end{center}
\begin{center}
Fig. 14
\end{center}
\end{figure}

The first case
($n=1$ and $l=2$)
for
an interior
vertex is impossible, because
the label of such a vertex, i.e., the product of labels of
corners, is
$a_m^{-1}p_2a_mp_1$ (if the vertex is a source) or
$b_0^{-1}p_2^\varphi b_0p_1^\varphi$ (if the vertex is a sink), where $p_1$ and
$p_2$ lie in $P$ and are not 1 (since the diagram is reduced)
and, therefore, the label of the vertex is not 1 by Corollary 2.1.;
thus this vertex can not be interior. The second and third
cases
($n=0$ and $l\in\{0,1\}$)
for
an interior
vertex are impossible by nearly the same reason:
they would imply an equality of the form $a_i^{\pm1}p_1=1$,
$b_i^{\pm1}p_1^\varphi=1$,
$p_2=1$, or $p_2^\varphi=1$, where $p_1\in P\ni p_2\ne1$.

Hence the curvature of any interior vertex or face is
nonpositive. The curvatures of exterior vertices are at
most two (this follows from formula (2) and the inequality $l\ge 2n$).
On the other hand, the total curvature must be four according to the weight test.
Therefore, the curvature of each exterior vertex is 2. This is possible only if
$n=0$, $l=0$, $p=0$, i. e. when all the corners at exterior vertices are the corners
of nonspecial digons. It follows that the diagram is a digon
with associated edges.
\\

{\large\textbf{6 The end of the proof.}}
\\

\sl
The aim of this section is to obtain the description of the element that conjugates
two elements from the subgroup $H$ in the group $\tilde{H}$ (presentation (3)). After that we will finish the proof.
 \\

\rm
\textbf{Theorem 6.1.} Suppose that elements $h$ and $g$ from the subgroup $H$ are conjugate in the group $\tilde{H}$, i. e.
$u^{-1}hu=g$ for some $u \in \tilde{H}$. Then such an irreducible form of
the element $u$ in the group $H\ast\langle t\rangle_\infty$ exists that it is a consecutive product
of elements from $H$ and words $u_i$, where $u_i$ is a label of a non-self-crossing path connecting exterior vertices
on the $i$-th spherical diagram from Theorem 4.1. (the number of all diagrams is $n$).

\textit{Proof.} The case $u \in H$ as a word from $H\ast
\langle t\rangle_\infty$ is obvious. Suppose the converse case.

Consider the process of the proof of  Theorem 4.1. for the equity $u^{-1}hug^{-1}=1$. We may assume that the word $u$ is
written as a product of the form $u_1h_1u_2h_2...u_kh_k$, where $u_i$ is
the word corresponding to the conjugation on the $i$-th sphere, i. e. $u_i$ is a label of a non-self-crossing path
connecting exterior vertices for all $i \in \{1,2,...k\}$. It is clear that the base of induction (the case $k=1$)
is proved after the <<thickening>> of paths, when we <<glue>> two paths with labels $u$ and $u^{-1}$.

Consider the diagram with the path corresponding to the word $u_i$ for $i \in \{1,2,...k\}$. In the process described above three
cases are possible.

\begin{enumerate}

\item A reducible pair is deleted.

\item An identity edge with different end points is pulled.

\item An identity edge with coinciding end points is pulled.

\end{enumerate}

Suppose that there is a subpath in the path with label $u_i$ relative to which the labels of faces of a certain reducible pair are inverse.
The considered subpath can be replaced by the path that runs the boundary of the right face (with respect to the
direction of the path with label $u_i$, Fig. 15).

\begin{figure}[h]
\begin{center}
\includegraphics[width=0.6\linewidth]{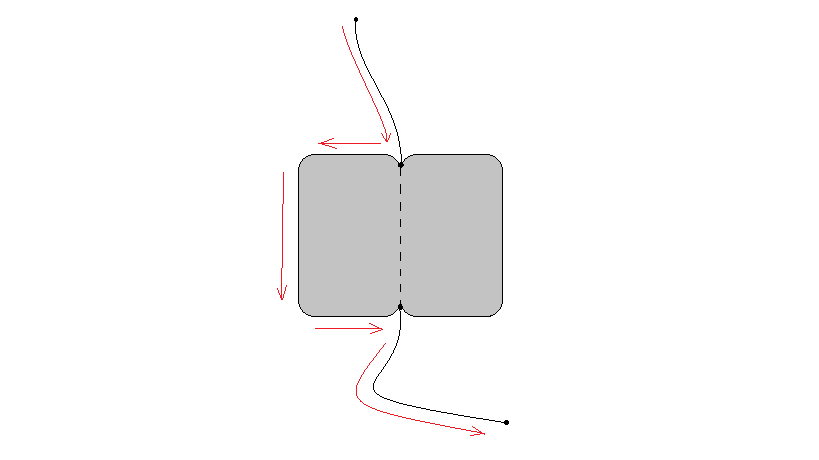}
\end{center}
\begin{center}
Fig. 15
\end{center}
\end{figure}

Obviously, after this replacement the label of the path does not change or it is
multiplied by an element $f_1\in H$ (in case, when the removed subpath is the beginning or the end
of the path) because labels of interior faces are equal to the identities in the group $\tilde{H}$.
Note that after the first operation self-crossings of the initial path may appear.

Consider the second case. If the identity edge connects two interior vertices then, obviously, the label of the path
does not change after the pulling, but a self-crossing of the path may appear. So we have two possibilities: two end points of
the identity edge are exterior vertices (Fig. 16a)) or only one of the end points is an exterior vertex
(Fig. 16b)).

\begin{figure}[h]
\begin{center}
\includegraphics[width=0.6\linewidth]{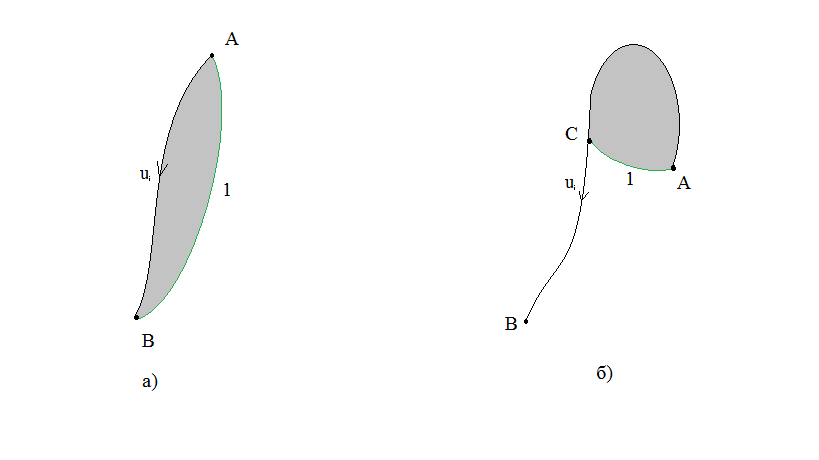}
\end{center}
\begin{center}
Fig. 16
\end{center}
\end{figure}

In case a) the vertices $A$ and $B$ are exterior. From Remark 4.2. it follows that the label of the cycle is equal
to the identity in the group $\tilde{H}$, hence $u_i=f_2, f_2\in H$.

In case b) only the vertex $A$ is exterior. In the same way as in case a) the label of the subpath $AC$ of the path $AB$ is an element
from $H$, hence after the pulling it is sufficient to consider the path $CB$. Let its label be $m$, then the equality $u_i=f_3m$ holds
in the group $\tilde{H}$ for some $f_3$ from $H$.

Our next aim is to remove self-crossings of the path.

\textit{Definition 6.1.} A subpath of the path with label $u_i$ is called \textit{a loop} if it is a cycle in the 1-skeleton of the diagram.

There are 4 types of self-crossings (a type
depends on the location of the exterior vertices with respect to the loop, Fig. 17 a), b), c), d)).

\begin{figure}[h]
\begin{center}
\includegraphics[width=0.6\linewidth]{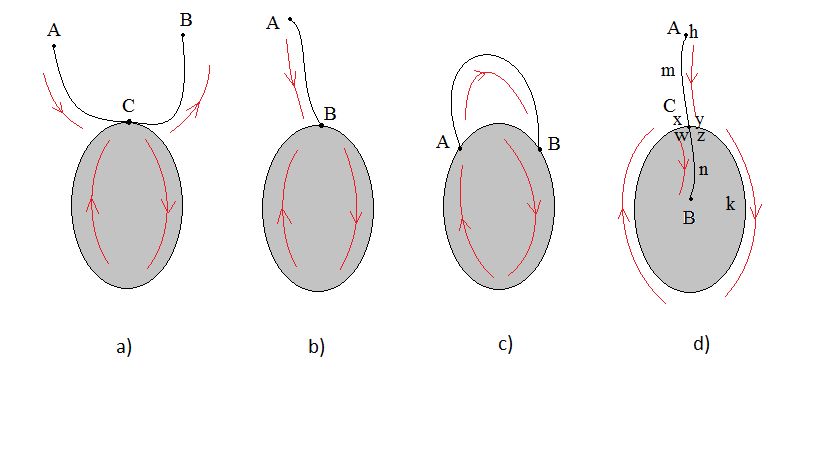}
\end{center}
\begin{center}
Fig. 17
\end{center}
\end{figure}

Every loop bounds two areas. In cases a), b), c) the exterior vertices are in the same area (may be, exactly on the loop).

In case a) Remark 4.2. implies that the label of the loop is the identity element in the group $H$, therefore the loop can be removed
from the path and the label of the path will not change.

In cases b) and c) the label of the loop equals to the identity in $H$ as above. Therefore, after removing the loop from
the path the number of self-crossings decreases and the label of the new path differs from the initial one by a multiplier
(left or right) $f_4$ from $H$. Indeed, the label of the path begins and ends with a label of an edge, hence after removing the loop the multiplier from $H$ may appear (this multiplier corresponds to the vertex $B$).

In case d) we can not use Remark 4.2. for the loop because each of the areas contains the exterior vertex. To consider this case we need some notation.
Let $m$ and $n$ be the labels of the subpaths $AC$ and $CB$ respectively, $k$ be the label of the loop $CC$, $x,y,z,w$
be the labels of four marked corners at the vertex $C$, $f_5$ be the label of the vertex $A$. Suppose that the subpath $AC$ does not
contain the exterior vertices (otherwise we can remove self-crossings of type b) as we did above). The label of the initial path equals
$mykxyzn$ and from Remark 4.2. it follows that the equality $mykxm^{-1}f_5=1$ holds in $\tilde{H}$. Hence the label of the
initial path is $f_5^{-1}myzn$, this differs from the label of the subpath $ACB$ (without the loop) by a multiplier from $H$.

Therefore after operations 1) and 2) from Theorem 4.1. we get two possibilities:

\begin{enumerate}

\item There is a path that connects exterior vertices, $u_i'$ is the label of this path, and the condition $u_i=h_1u_i'h_2$ holds in the group $\tilde{H}$, where $h_1$ and $h_2$ are elements of $H$.

\item The conjugation by the element $u_i$ is a conjugation by an element from $H$.

\end{enumerate}

Consider operation 3) (Fig. 18). Without loss of generality we may assume that the exterior vertices are contained in different areas.
Let $A$ and $B$ be the end points of the considered path (its label is $u_i$), $r_1$ and $r_2^{-1}$ be the labels of the respective vertices,
$C$ be the beginning and the end of the identity edge, $u_i'$ and $u_i''$ be the labels of the subpaths $AC$ and $CB$ of the path $AB$.
The reader will easily prove that the following equations hold in $\tilde{H}$:

$u_i'^{-1}r_1u_i'=z^{-1}x^{-1}$;

$u_i''^{-1}wyu_i''=r_2$;

$z^{-1}w^{-1}=xy$;

$(wz)(z^{-1}x^{-1})(xy)=wy$

Therefore $u_i$ can be replaced by the product $u_i'z^{-1}w^{-1}u_i''$, where $u_i'$ and $u_i''$ correspond to two new spherical diagrams.
This completes the proof.

\begin{figure}[h]
\begin{center}
\includegraphics[width=0.6\linewidth]{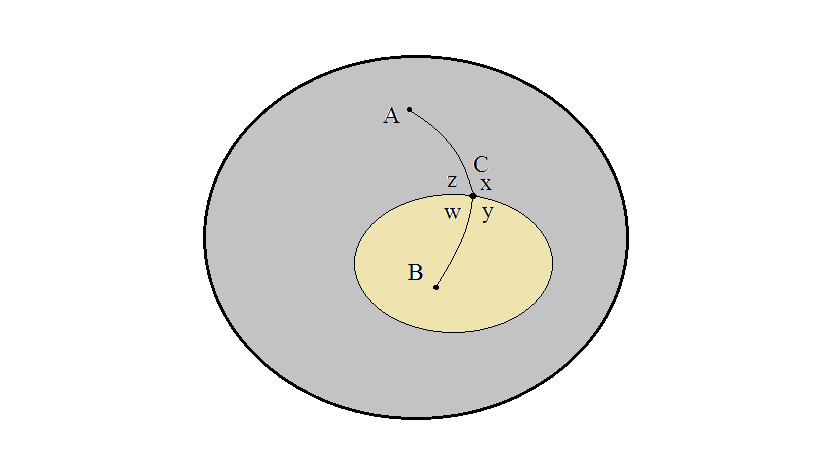}
\end{center}
\begin{center}
Fig. 18
\end{center}
\end{figure}

\textit{Definition 6.2.} Consider a word $v$ from the group $H\ast\langle t\rangle_\infty$. \textit{A prefix} of this word is an element
of  $H\ast\langle t\rangle_\infty$, such that it is a product of some first letters of the word $v$.

\textbf{Corollary 6.1.} If presentation (3) has the particular type of presentation (1) and elements $g$ and $h$ of the group $H$ are conjugate
in $\tilde{H}$ by an element $u$, then there exists a presentation of $u$ as a word from $H\ast\langle t\rangle_\infty$, that for every prefix $v$
of the word $u$ the condition $v^{-1}hv\in H$ holds in $\tilde{H}$.

\textit{Proof.} Indeed, Theorem 5.1. claims that the only $\varphi$-reduced Howie diagram without exterior faces with two exterior vertices over presentation (1) is the diagram of one degenerated
digon. Hence, every Howie diagram from Theorem 4.1. corresponds to the conjugation of elements of the subgroups $P$ and $P^{\varphi}$ by the element $t$. This assertion and Theorem 6.1. prove Corollary 6.1.

Let us finish the proof of Theorem 1.1. Note that if $w=gt$, where $g$ is an element of the group $G$, then we have an isomorphism:

$\tilde{G}=\langle G,t|(gt)^k =1\rangle\cong_{\phi} G\ast\langle x
\rangle_{k}$

The structure of the isomorphism $\phi$ is:

$\phi(g)=g,\forall g \in G$, $\phi(t)=g^{-1}x$

Therefore the malnormality of $G$ in $\tilde{G}$ is obvious.

Now suppose that there are two or more letters $t$ in the word $w$. We apply Lemma 2.1. and all the statements from Sections 4 and 5.

Let the elements $h^{(0)}$ and $g^{(0)}$ from the (0)-multiplier of the free product $H$ be the isomorphic images (the natural isomorphism between
presentations (0) and (1)) of the elements $h$ and $g$ from $G$ and the equality $u^{-1}h^{(0)}u=g^{(0)}$ holds in $\tilde{H}$ for some word
$u \in H\ast\langle t\rangle_\infty$, where $u$ is an irreducible word. From Corollary 6.1. we can suppose that the element $h^{(0)}$ conjugated by every prefix
of the word $u$ is an element of the subgroup $H$ and all the conjugations by the element $t^{\pm 1}$ correspond to the diagrams of one degenerated
digon. Let $u$ be $a_1t^{k_1}...a_nt^{k_n}a_{n+1}$, where $a_i \in H$, $k_i$ are integers $\forall i \in \{1,...,n+1\}$. It is clear that $\sum k_i=0$ because
the result of all the conjugations is the element of the (0)-multiplier. Now let us prove by induction over $\sum |k_i|$ that $u \in H$.

The base $\sum |k_i|=0$ is obvious.

Consider the inductive step. Without loss of generality we can assume that there is a subword $t^{-1}a_it$ for some $i$ in the word $u$ because $\sum k_i=0$.
Then for some $h \in H$ the following equalities hold in $\tilde{H}$:

$tht^{-1}=h_1$, where $h_1 \in H$;

$a_i^{-1}h_1a_i=h_2$, where $h_2 \in H$;

$t^{-1}h_2t=h_3$, where $h_3 \in H$

Obviously, there are no multipliers from $G^{(s)}$ in the irreducible form of the element $h_1$, hence there are no such multipliers in the irreducible
form of the element $a_i$, otherwise the last equality is false. Therefore $t^{-1}a_it \in H$. We decreased $\sum |k_i|$ and
the step is proven.

Since $H$ is the free product of several copies of the group $G$, we obtain that $u \in H$ implies $u \in G^{(0)}$. This completes the proof.
\\

{\large\textbf{7 The second proof.}}
\\

Let us begin with D. Osin's result that was mentioned in the introduction.

\textbf{Theorem 7.1.} \cite{Os} Let $M$ be a group, $\{H_\lambda\}_{\lambda\in\Lambda}$ is a collection of subgroups of the group $M$.
Suppose that $M$ is relatively finitely presented with respect to $\{H_\lambda\}_{\lambda\in\Lambda}$, and the relative Dehn function of the group $M$
with respect to $\{H_\lambda\}_{\lambda\in\Lambda}$ is well-defined, i. e. for all natural $n$ a finite value of Dehn function  $f_D(n)$ exists.
Then the following statements are true:

1) for all $m_1, m_2 \in M$ the intersection $H_\lambda^{m_1}\cap H_\mu^{m_2}$ is finite if $\lambda \neq \mu$.

2) The intersection $H_\lambda^m \cap H_\lambda$ is finite for all $m \notin H_\lambda$.

Then consider the following theorem.

\textbf{Theorem 7.2.} \cite{Kl-Lu} If a word $w\in G\ast \langle t\rangle_\infty$ is unimodular and $k\geq2$, then the group $\tilde{G}$ defined
by relative presentation (0) contains $G$ as a subgroup (naturally embedded). Moreover, if the group $G$ is involution free or $k\geq3$, then
$\tilde{G}$ is respectively hyperbolic with respect to $G$.

Thereby under the conditions of Theorem 7.2. presentation (0) is aspherical, and the relative hyperbolicity (i.e. the linearity
of the respective Dehn function) implies that we can apply Theorem 7.1.

Consider it in detail. Let $u\in \tilde{G}$ be an element that is not contained in the subgroup $G$. Then the intersection $G\cap G^u$ is finite.
Suppose that $g=h^u$, where $g$ and $h$ are nontrivial elements of finite order $s$ from $G\cap G^u$. From Theorem 4.1. it follows that there are $n \in \mathbb{N}$
reduced (or $\varphi$-reduced, this is the same in case of presentation (0)) spherical Howie diagrams over presentation (0) without exterior faces with two exterior vertices with
labels $g_i'$ and $g_{i+1}^{-1}$ from $G$, $i \in \{0,...n-1\}$, where elements $g_l$ and $g_l'$ are conjugate in $G$ for all $l \in \{1,...n-1\}$, and $g'_0=h, g_n=g$.
Obviously, the labels of the exterior vertices have order $s$. Let $u_i, i \in \{0,...n-1\}$ be a label of some non-self-crossing path (from Theorem 6.1.),
that connects the exterior vertices on the $i$-th diagram. Theorem 6.1. implies that there is such an $i, i \in \{0,...n-1\}$, that $u_i\in \tilde{G}\backslash G$
(otherwise, $u$ is from subgroup $G$).

Let us cut the path with label $u_i$ on the $i$-th diagram and take $s$ copies of such a disk diagram. Let us number the copies from $0$ to $s-1$ and <<glue>> the path
with label $u_i^{-1}$ of the $j$-th copy with the path with label $u_i$ of the (j+1)-th copy (subscripts are modulo s) such a way that the edges with labels
$g_i'$ and $g_{i+1}^{-1}$ become common for <<glued>> copies. Thereby we obtain a reduced spherical diagram without exterior vertices and faces because the initial spherical diagram is reduced
and the orders of the elements $g_i'$ and $g_{i+1}^{-1}$ are equal to $s$. The existence of such a diagram contradicts to the asphericity of presentation (0), hence the intersection
$G\cap G^u$ is trivial.

The following theorem is also very close to the solved problem.

\textbf{Theorem 7.3.} \cite{FR-05} Let the relative 2-complex $(L,K)$ be relatively aspherical, i. e. the mapping $\pi_2(K\cup L^{(1)},K)\rightarrow\pi_2(L,K)$ is injective. Then every
finite subgroup $\pi_1(L)$ is contained only in one group conjugated to $\pi_1(K)$.

The main idea is that the asphericity in terms of complexes follows from the asphericity in terms of Howie diagrams, but the problem is that we can not formally use Theorem 7.3. because
it is not obvious that the subgroups $G$ and $G^u$ do not coincide.

\end{document}